\documentclass[preprint,12pt]{elsarticle}

\usepackage{amssymb}
\usepackage{amsthm}
\usepackage{graphics}
\usepackage{epsfig}
\usepackage{graphicx}
\usepackage{color}
\usepackage{hyperref}
\usepackage{a4wide}
\usepackage{amsmath}
\usepackage{amsfonts}
\usepackage{enumerate}
\newcommand{\pf}{\noindent {\bf Proof: }}

\newtheorem{lemma}{\bf Lemma}[section]
 
\newtheorem{theorem}{\bf Theorem}[section]
\newtheorem{proposition}[lemma]{\bf Proposition}
\newtheorem{corollary}[lemma]{\bf Corollary}

\newtheorem{remark}{\bf Remark}[section]

\journal{arxiv.org}

\begin{document}

\begin{frontmatter}



\title{Partial Domination in Graphs}



\author{Angsuman Das\corref{cor1}}
\ead{angsumandas@sxccal.edu}

\address{Department of Mathematics,\\ St. Xavier's College, Kolkata, India.\\angsumandas@sxccal.edu}
\cortext[cor1]{Corresponding author}


\begin{abstract}
Let $G=(V,E)$ be a graph. For some $\alpha$ with $0<\alpha \leq 1$, a subset $S$ of $V$ is said to be a $\alpha$-partial dominating set if $|N[S]|\geq \alpha |V|$. The size of a smallest such $S$ is called the $\alpha$-partial domination number and is denoted by $\mathsf{pd}_\alpha(G)$. In this paper, we introduce $\alpha$-partial domination number in a graph $G$ and study different bounds on the partial domination number of a graph $G$ with respect to its order, maximum degree, domination number etc.,  Moreover, $\alpha$-partial domination spectrum is introduced and Nordhaus-Gaddum bounds on the partial domination number are studied. 
\end{abstract}

\begin{keyword}
$\alpha$-domination \sep connected graph \sep edge-critical graphs
\MSC[2008] 05C69

\end{keyword}

\end{frontmatter}


\section{Introduction}
The domination in graphs has been an active area of research from the time of its inception. Two domination books \cite{domination-book,advanced-domination-book} provide a comprehensive report of vastness of the area of the domination and its relation to other graph parameters. Many variations e.g., \cite{fair-dom,angsu-coeff-dom,angsu-icmc-17,alpha-dom} etc.,  of the domination problem can be found in literature most of which are motivated by many real-life scenarios. 

Consider the following scenario.   Imagine that you are the curator of an art museum and you wish to determine the minimum number of guards you need to guard the exhibits.  A guard can guard an exhibit that he/she is standing near, and any exhibit in the museum that they can clearly see.  In order to model the security situation, you would construct a graph $G$ as following:  Each vertex represents an exhibit location and two vertices $u$ and $v$ are adjacent if and only if  the locations they represent are visible from each other, that is, a person standing at the exhibit modelled by vertex $u$ can clearly see the location of the exhibit modelled by vertex $v$, and vice versa.   Suppose that security requirements mandate that a staff of guards are positioned at locations such that every art exhibit is protected by a guard that can see it, and budget restrictions make it desirable to hire as few guards as possible.  In this case, the most economical solutions, that is, the minimum guards for possible guard location configurations, correspond to the  $\gamma$-sets.  Suppose that due to budgetary concerns, as  curator, you are  strictly limited to hiring exactly $\gamma$ guards.

While this is the optimum solution economically,  in a practical sense it leaves much to be desired.  There will be days when guards are ill, guards need a day off or some of them institute a labour action and go on strike.  As curator you can now at most secure a fraction or part of the exhibits and keep the rooms containing unguarded exhibits locked for that day. It is with this problem in mind, that we introduce in this paper the concept of the  {\it partial domination} in a graph. A closely related problem based on algorithmic viewpoint can be found in \cite{t-domination}.

Let $G=(V,E)$ be a graph. For some $\alpha$ with $0<\alpha \leq 1$, a subset $S$ of $V$ is said to be a $\alpha$-partial dominating set if $|N[S]|\geq \alpha |V|$. The size of a smallest such $S$ is called the $\alpha$-partial domination number and is denoted by $\mathsf{pd}_\alpha(G)$. Clearly $1 \leq \mathsf{pd}_\alpha(G)\leq \gamma(G)$ and $\mathsf{pd}_{1}(G)=\gamma(G)$. Also, $\alpha_1 < \alpha_2$ implies $\mathsf{pd}_{\alpha_1}(G)\leq \mathsf{pd}_{\alpha_2}(G)$.

\section{Some Basic Results}
We start with some basic results. As they are straightforward, they are given either without proof or with a minimalistic proof.
{\proposition \label{basic-bound-1} Let $G$ be a graph on $n$ vertices. Then $\mathsf{pd}_\alpha(G)=1$ for all $\alpha \in (0, \frac{\Delta+1}{n}]$.}

{\proposition \label{basic-bound-2} Let $G$ be a graph on $n$ vertices. Then $\mathsf{pd}_\alpha(G)=\gamma(G)$ for all $\alpha \in (1-\frac{1}{n}, 1]$.}

{\proposition \label{paths-and-cycles} $\mathsf{pd}_\alpha(C_n)=\mathsf{pd}_\alpha(P_n)=\lceil \frac{n\alpha}{3}\rceil$}\\
\pf Let $S$ be a $\mathsf{pd}_\alpha$-set of $C_n$. Then $|N[S]|\geq \lceil n\alpha\rceil$. To dominate $\lceil n\alpha\rceil$ vertices in $C_n$, we need at least $\left\lceil \frac{\lceil n\alpha\rceil}{3}\right\rceil$ vertices. Thus $\mathsf{pd}_\alpha(C_n)=|S|=\left\lceil \frac{\lceil n\alpha\rceil}{3}\right\rceil=\lceil \frac{n\alpha}{3}\rceil$. Similarly,
$\mathsf{pd}_\alpha(P_n)=\lceil \frac{n\alpha}{3}\rceil$.\qed

{\proposition $\mathsf{pd}_\alpha(K_n)=1$ and for $m\geq n$, $\mathsf{pd}_\alpha(K_{m,n})=
\left\lbrace\begin{array}{ll}
1 & \mbox{ if }0<\alpha \leq\frac{m+1}{m+n},\\
2 & \mbox{ if }\alpha < \frac{m+1}{m+n}\leq 1].
\end{array}\right.
$}

\section{Some Bounds on Partial Domination Number}
In this section, we study some bounds on partial domination number of  a graph with respect to other graph parameters.
{\proposition \label{=1} Let $G=(V,E)$ be a graph on $n$ vertices. Then $\mathsf{pd}_\alpha(G)=1$ if and only if there exists $v \in V$ such that $deg(v)\geq \lceil n \alpha \rceil -1$.}\\
\pf $\mathsf{pd}_\alpha(G)=1$ if and only if there exists a vertex $v \in V$ such that $|N[v]|\geq n\alpha$, i.e., $deg(v)\geq \lceil n \alpha \rceil -1$.\qed

{\proposition \label{lower-bound-1} Let $G=(V,E)$ be a graph of order $n$ and maximum degree $\Delta$ such that $\Delta < \lceil n \alpha \rceil -1$. Then $\frac{n \alpha}{\Delta +1} \leq \mathsf{pd}_\alpha(G)\leq \lceil n \alpha \rceil -\Delta$.}\\
\pf Let $S$ be a $\mathsf{pd}_\alpha$-set in $G$. Then $$n \alpha \leq |N[S]| \leq \sum_{v \in S} deg(v) + |S|\leq (\Delta +1)|S|=(\Delta +1)\mathsf{pd}_\alpha$$ and hence the lower bound  follows.

For the upper bound, let $v$ be a vertex of maximum degree in $G$. Then $v$ dominates $\Delta +1$ vertices. Then $v$ along with other $\lceil n\alpha \rceil -(\Delta+1)$ vertices outside $N[v]$ forms a $\alpha$-partial dominating set of $G$. Thus $\mathsf{pd}_\alpha(G)\leq 1+\lceil n \alpha \rceil -(\Delta+1)=\lceil n \alpha \rceil -\Delta$.\qed

{\proposition \label{upper-bound-1} Let $G$ be a graph with domination number $\gamma$. Then $\mathsf{pd}_\alpha(G)\leq \left\lceil \dfrac{\gamma}{\lfloor \frac{1}{\alpha} \rfloor} \right\rceil$.}\\
\pf Let $D$ be a $\gamma$-set of $G$ and set $t=\lfloor \frac{1}{\alpha}\rfloor$. Let $D_1,D_2,\ldots, D_t$ be a partition of $D$ such that $|D_i|\leq \left\lceil \dfrac{\gamma}{\lfloor \frac{1}{\alpha} \rfloor} \right\rceil$ for all $i$. Thus we have $N[D]=N[D_1]\cup N[D_2]\cup \cdots \cup N[D_t]$. Let $n$ be the order of $G$. Then 
$$n=|N[D]|\leq \sum_{i=1}^t |N[D_i]|\leq t|N[D_j]|, \mbox{ where }|N[D_j]|=\max_i |N[D_i]|$$
$$i.e.,~|N[D_j]|\geq \dfrac{n}{t}=\dfrac{n}{\lfloor \frac{1}{\alpha}\rfloor}\geq n\alpha$$
Thus, $N[D_j]$ is an $\alpha$-partial dominating set of $G$ and hence $\mathsf{pd}_\alpha(G)\leq \left\lceil \dfrac{\gamma}{\lfloor \frac{1}{\alpha} \rfloor} \right\rceil$.\qed

{\corollary If $G$ is a graph of order $n$ without any isolated vertex, then $\mathsf{pd}_\alpha(G)\leq \left\lceil \dfrac{n}{2\lfloor \frac{1}{\alpha} \rfloor} \right\rceil$.}\\
\pf Since $G$ does not have any isolated vertex, $\gamma \leq n/2$. Thus the corollary follows from the previous theorem.\qed

{\corollary \label{1/gamma} If $G$ is a graph with domination number $\gamma$ and $\alpha \in (0,1/\gamma]$, then $\mathsf{pd}_\alpha(G)=1$.}\\
\pf Since $\alpha \leq 1/\gamma$, we have $\gamma \leq 1/\alpha$, i.e., $\gamma \leq \lfloor 1/\alpha \rfloor$. Hence $\dfrac{\gamma}{\lfloor 1/\alpha \rfloor}\leq 1$, i.e., $\left\lceil \dfrac{\gamma}{\lfloor \frac{1}{\alpha} \rfloor} \right\rceil=1$ and thus by Theorem \ref{upper-bound-1}, we have $\mathsf{pd}_\alpha(G)=1$.\qed

{\theorem \label{alpha-1-alpha} Let $G$ be a graph with domination number $\gamma(G)$. Then for all $\alpha \in (0,1)$, $\mathsf{pd}_\alpha(G)+\mathsf{pd}_{1-\alpha}(G)\leq \gamma + 1$.}\\
\pf Let $S$ be a $\gamma(G)$-set and $\alpha \in (0,1)$. Let $S_1$ be a subset of $S$ with $|N[S_1]|\geq n\alpha$ such that $S_1$ is a minimal subset of $S$ with this property. Clearly $\mathsf{pd}_\alpha(G)\leq |S_1|$. Let $S_2=S\setminus S_1$ and $v \in S_1$. Since $S_1$ is minimal with respect to the above property, we have $|N[S_1\setminus \{v\}]|<n\alpha$. Now, as $S=\left(S_1\setminus \{v\}\right)\cup \left( S_2\cup \{v\}\right)$, we get $$n=|V|=|N[S]| \leq |N[\left(S_1\setminus \{v\}\right)]|+|N[\left( S_2\cup \{v\}\right)]|<n\alpha +|N[\left( S_2\cup \{v\}\right)]|$$
$$\mbox{i.e., }|N[\left( S_2\cup \{v\}\right)]|>n(1-\alpha)$$
Thus $S_2\cup \{v\}$ is an $(1-\alpha)$-partial dominating set of $G$ and  $\mathsf{pd}_{1-\alpha}(G)\leq |S_2\cup \{v\}|=|S_2|+1$. Hence, $$\mathsf{pd}_\alpha(G)+\mathsf{pd}_{1-\alpha}(G)\leq |S_1|+|S_2|+1=|S|+1=\gamma+1.$$\qed

In fact, it is possible to find a generalization of Theorem \ref{alpha-1-alpha} in a natural way.
{\theorem \label{alpha1+alpha2+alphak} Let $G$ be a graph with domination number $\gamma$. For any positive integer $k \geq 2$, with $\alpha_1+\alpha_2+\cdots + \alpha_k \leq 1$ and $\alpha_i \in (0,1)$ for all $i$, $\mathsf{pd}_{\alpha_1}(G)+\mathsf{pd}_{\alpha_2}(G)+\cdots +\mathsf{pd}_{\alpha_k}\leq \frac{k}{2}(\gamma + 1)$.}\\
\pf We prove it by induction on $k$. For $k=2$, $\alpha_1+\alpha_2 \leq 1$. Hence, by Theorem \ref{alpha-1-alpha}, $\mathsf{pd}_{\alpha_1}(G)+\mathsf{pd}_{\alpha_2}(G) \leq \mathsf{pd}_{\alpha_1}(G)+\mathsf{pd}_{1-\alpha_1}(G)\leq \gamma + 1$. Assume that $k >2$ and the theorem holds for integers less than $k$. Then at least one value of $\alpha_i$ must satisfy $\alpha_i \leq \frac{1}{2}$. Without loss of generality, let $\alpha_k \leq \frac{1}{2}$. Therefore, by Theorem \ref{upper-bound-1}, $\mathsf{pd}_{\alpha_k}(G) \leq \left\lceil \dfrac{\gamma}{2} \right\rceil$. Finally, using the induction hypothesis, we get $$[\mathsf{pd}_{\alpha_1}(G)+\mathsf{pd}_{\alpha_2}(G)+\cdots +\mathsf{pd}_{\alpha_{k-1}}(G)]+\mathsf{pd}_{\alpha_k}(G) \leq \frac{(k-1)}{2}(\gamma + 1)+\left\lceil \dfrac{\gamma}{2} \right\rceil$$

$$\mbox{i.e., }\mathsf{pd}_{\alpha_1}(G)+\mathsf{pd}_{\alpha_2}(G)+\cdots +\mathsf{pd}_{\alpha_k}(G) \leq \frac{(k-1)}{2}(\gamma + 1)+ \dfrac{\gamma}{2}+\dfrac{1}{2}=\frac{k}{2}(\gamma + 1).$$\qed

{\theorem Let $G$ be a graph with components $G_1,G_2,\ldots,G_k$. Then $$ \mathsf{pd}_\alpha(G)\leq \sum_{i=1}^k \mathsf{pd}_\alpha(G_i)$$. }\\
\pf Let $S_i$ be a $\mathsf{pd}_\alpha(G_i)$-set of $G_i$, for $i=1,2,\ldots,k$. Then $|N[S_i]|\geq \alpha |V(G_i)|$, for $i=1,2,\ldots,k$. Let $S=S_1\cup S_2 \cup \cdots \cup S_k$. Thus $$|N[S]|=\sum_{i=1}^k |N[S_i]|\geq \alpha \sum_{i=1}^k |V(G_i)|=\alpha |V(G)|,$$
and $S$ is a $\alpha$-partial dominating set of $G$ and hence, $$\mathsf{pd}_\alpha(G)\leq |S|=\sum_{i=1}^k |S_i|=\sum_{i=1}^k \mathsf{pd}_\alpha(G_i).$$ \qed

\section{Vertex and Edge Removal}
In this section , we focus on effect of removal and addition of edges and vertices of a graph on its partial domination number.
{\theorem Let $G=(V,E)$ be a graph and $e \in E$. Then $\mathsf{pd}_\alpha(G) \leq \mathsf{pd}_\alpha(G-e)\leq \mathsf{pd}_\alpha(G)+1$.}\\
\pf Clearly $\mathsf{pd}_\alpha(G) \leq \mathsf{pd}_\alpha(G-e)$. Thus we prove the other part of the inequality. Let $S$ be a $\mathsf{pd}_\alpha(G)$-set and $e=xy$ where $x,y \in V$. If $x,y \in S$, then $N_G[S]=N_{G-e}[S]$ and hence $S$ is an $\alpha$-partial dominating set of $G-e$, i.e., $\mathsf{pd}_\alpha(G-e) \leq \mathsf{pd}_\alpha(G)$. Similarly, if $x \not\in S$ and $y \not\in S$, then $N_G[S]=N_{G-e}[S]$ and hence $\mathsf{pd}_\alpha(G-e) \leq \mathsf{pd}_\alpha(G)$. Finally, if $x \not\in S$ and $y \in S$, then $S\cup \{x\}$ is an $\alpha$-partial dominating set of $G-e$, i.e., $\mathsf{pd}_\alpha(G-e) \leq \mathsf{pd}_\alpha(G)+1$. Combining all the cases, we get the upper bound. \qed

{\theorem Let $G=(V,E)$ be a graph and $v \in V$. Then $\mathsf{pd}_\alpha(G)-1 \leq \mathsf{pd}_\alpha(G-v)\leq \mathsf{pd}_\alpha(G)+deg_G(v)-1$.}\\
\pf Let $S$ be a $\mathsf{pd}_\alpha(G)$-set. If $v \not\in N[S]$, then $S$ is an $\alpha$-partial dominating set of $G-v$, as $|N_{G-v}[S]|=|N_G[S]|\geq |V|\alpha > (|V|-1)\alpha$ and hence $\mathsf{pd}_\alpha(G-v)\leq \mathsf{pd}_\alpha(G)$. 

If $v \in N[S]\setminus S$, then $N_{G-v}[S]=N_G[S]\setminus \{v\}$, i.e., $|N_{G-v}[S]|=|N_G[S]|-1\geq |V|\alpha -1$. Let $u \in V \setminus N_{G-v}[S]$ such that $u \neq v$ and set $S'=S\cup \{u\}$. Then
$$|N_{G-v}[S']|\geq |N_{G-v}[S]|+1 \geq |V|\alpha \geq (|V|-1)\alpha$$
i.e., $S'$ is an $\alpha$-partial dominating set of $G-v$ and hence $\mathsf{pd}_\alpha(G-v)\leq |S'|=\mathsf{pd}_\alpha(G)+1$.

If $v \in S$, then $N_{G-v}[S\setminus \{v\}]\supseteq N_G[S]\setminus N[v]$, i.e., $$|N_{G-v}[S\setminus \{v\}]|\geq |N_G[S]|-(deg_G(v)+1)\geq |V|\alpha -deg_G(v)-1.$$
Now let $T$ be a collection of $deg_G(v)$ many vertices in $(V\setminus \{v\}) \setminus N_{G-v}[S\setminus \{v\}]$ and let $S_1=(S\setminus \{v\})\cup T$. Then $N_{G-v}[S_1]=N_{G-v}[S\setminus \{v\}]\cup N_{G-v}[T]$ and hence $$|N_{G-v}[S_1]| \geq |N_{G-v}[S\setminus \{v\}]|+|T|\geq (|V|\alpha -deg_G(v)-1) + deg_G(v)=|V|\alpha -1 \geq (|V|-1)\alpha.$$
Thus $S_1$ is an $\alpha$-partial dominating set of $G-v$, i.e., $$\mathsf{pd}_\alpha(G-v)\leq |S_1|=|S|-1+deg_G(v)=\mathsf{pd}_\alpha(G)+deg_G(v)-1.$$
Combining all the above three cases, we have the proposed upper bound.

For the lower bound, let $S$ be a $\mathsf{pd}_\alpha(G-v)$-set. Then $|N_{G-v}[S]|\geq (|V|-1)\alpha$. Let $u \in V\setminus N_{G-v}[S]$ and set $S_1=S\cup \{u\}$. Then $$|N_G[S_1]|\geq |N_{G-v}[S]|+1\geq (|V|-1)\alpha +1\geq |V|\alpha.$$
Thus, $S_1$ is an $\alpha$-partial dominating set of $G$ and hence, $\mathsf{pd}_\alpha(G)\leq |S_1|=|S|+1=\mathsf{pd}_\alpha(G-v)+1$, i.e., $\mathsf{pd}_\alpha(G)-1 \leq \mathsf{pd}_\alpha(G-v)$.\qed

We call a graph $G$, $\alpha$-partial domination vertex critical, or just $\mathsf{pd}_\alpha$-vertex critical if for any $v \in V$, $\mathsf{pd}_\alpha(G-v)<\mathsf{pd}_\alpha(G)$. In the light of the above theorem, if a graph $G$ is $\mathsf{pd}_\alpha$-vertex critical, then  $\mathsf{pd}_\alpha(G-v)=\mathsf{pd}_\alpha(G)-1$. 

{\theorem If $G=(V,E)$ be a $\mathsf{pd}_\alpha$-vertex critical graph, then for every vertex $v \in V$, there exists a $\mathsf{pd}_\alpha(G)$-set $S$ containing $v$ such that $pn_G[v,S]=\{v\}$.}\\
\pf Let $T$ be a $\mathsf{pd}_\alpha(G-v)$-set. Then $|N_{G-v}[T]|\geq (|V|-1)\alpha$. Thus $$|N_G[T\cup \{v\}]|\geq (|V|-1)\alpha +1\geq |V|\alpha,$$ and hence $T\cup \{v\}$ is an $\alpha$-partial dominating set of $G$, i.e., $\mathsf{pd}_\alpha(G)\leq |T\cup \{v\}|=\mathsf{pd}_\alpha(G-v)+1$. Moreover as $G$ is $\mathsf{pd}_\alpha$-vertex critical graph, we have $\mathsf{pd}_\alpha(G-v)=\mathsf{pd}_\alpha(G)-1$ for all $v \in V$. Thus $\mathsf{pd}_\alpha(G)= |T\cup \{v\}|$, i.e., $S=T \cup \{v\}$ is a $\mathsf{pd}_\alpha(G)$-set containing $v$.

If $T \cap N(v) \neq \emptyset$, then $|N_G[T]|=|N_{G-v}[T]|+1\geq (|V|-1)\alpha +1\geq |V|\alpha$, i.e., $T$ is an $\alpha$-partial dominating set of $G$. This contradicts $\mathsf{pd}_\alpha(G-v)<\mathsf{pd}_\alpha(G)$ and hence $T \cap N(v) = \emptyset$. Thus $pn_G[v,S]=\{v\}$. \qed

\section{Nordhaus-Gaddum Bounds}
In this section, we study some Nordhaus-Gaddum bounds on partial domination number of a graph $G$. We start with recalling some known Nordhaus-Gaddum bounds on domination number of a graph $G$.
{\proposition [Cockayne and Hedeitniemi] \label{cockayne-hedeitniemi} For any graph $G$, $\gamma(G)+\gamma \left(\overline{G}\right)\leq n+1$ with equality if and only if $G=K_n$ or $G=\overline{K_n}$.\qed }

{\proposition[Laskar and Peters] \label{laskar-peters} For connected graphs $G$ and $\overline{G}$, $\gamma(G)+\gamma \left(\overline{G}\right)\leq n$ with equality if and only if $G=P_4$.\qed }

{\proposition[Joseph and Arumugam] \label{joseph-arumugam} For graphs $G$ and $\overline{G}$ without isolated vertices, $\gamma(G)+\gamma \left(\overline{G}\right)\leq \lfloor n/2 \rfloor +2$.\qed }

{\theorem \label{NG1} For any graph $G$, $\mathsf{pd}_\alpha(G)+\mathsf{pd}_\alpha(\overline{G})\leq \left\lceil \dfrac{n}{\lfloor \frac{1}{\alpha} \rfloor} \right\rceil +1$.}\\
\pf From Theorem \ref{upper-bound-1}, we have $\mathsf{pd}_\alpha(G)\leq \left\lceil \dfrac{\gamma(G)}{\lfloor \frac{1}{\alpha} \rfloor} \right\rceil$. Thus $$\mathsf{pd}_\alpha(G)+\mathsf{pd}_\alpha(\overline{G})\leq \left\lceil \dfrac{\gamma(G)}{\lfloor \frac{1}{\alpha} \rfloor} \right\rceil+ \left\lceil \dfrac{\gamma(\overline{G})}{\lfloor \frac{1}{\alpha} \rfloor} \right\rceil \leq \left\lceil \dfrac{\gamma(G) + \gamma(\overline{G})}{\lfloor \frac{1}{\alpha} \rfloor} \right\rceil+1.$$

From Proposition \ref{cockayne-hedeitniemi}, we get $\gamma(G)+\gamma(\overline{G})\leq n$ except when $G=K_n$ or $G=\overline{K_n}$. Thus apart from these two cases, we have $\mathsf{pd}_\alpha(G)+\mathsf{pd}_\alpha(\overline{G})\leq \left\lceil \dfrac{n}{\lfloor \frac{1}{\alpha} \rfloor} \right\rceil +1$. 

Now, consider the case when $G=K_n$ (or $\overline{K_n}$). Then $\mathsf{pd}_\alpha(G)=1$ and $\mathsf{pd}_\alpha(\overline{G})=\lceil n\alpha \rceil$. Thus, $\mathsf{pd}_\alpha(G)+\mathsf{pd}_\alpha(\overline{G})=\lceil n\alpha \rceil +1 \leq \left\lceil \dfrac{n}{\lfloor \frac{1}{\alpha} \rfloor} \right\rceil +1$. Combining all these cases, we have the theorem.\qed

{\theorem For connected graphs $G$ and $\overline{G}$, $\mathsf{pd}_\alpha(G)+\mathsf{pd}_\alpha(\overline{G})\leq \left\lceil \dfrac{n-1}{\lfloor \frac{1}{\alpha} \rfloor} \right\rceil +1$.}\\
\pf Similar to the proof of Theorem \ref{NG1}, we have $\mathsf{pd}_\alpha(G)+\mathsf{pd}_\alpha(\overline{G})\leq \left\lceil \dfrac{\gamma(G) + \gamma(\overline{G})}{\lfloor \frac{1}{\alpha} \rfloor} \right\rceil+1.$

From Proposition \ref{laskar-peters}, we get $\gamma(G)+\gamma \left(\overline{G}\right)\leq n-1$ except when $G=P_4$. Thus, apart from the case when $G=P_4$, we have $\mathsf{pd}_\alpha(G)+\mathsf{pd}_\alpha(\overline{G})\leq \left\lceil \dfrac{n-1}{\lfloor \frac{1}{\alpha} \rfloor} \right\rceil +1$.

Now, consider the case when $G=\overline{G}=P_4$. If $0<\alpha \leq 3/4$, then by Proposition \ref{paths-and-cycles}, $\mathsf{pd}_\alpha(G)=\mathsf{pd}_\alpha(\overline{G})=1$, i.e., $\mathsf{pd}_\alpha(G)+\mathsf{pd}_\alpha(\overline{G})=2\leq \left\lceil \dfrac{4-1}{\lfloor \frac{1}{\alpha} \rfloor} \right\rceil +1$. 

If $3/4 < \alpha \leq 1$, we have $\mathsf{pd}_\alpha(G)=\mathsf{pd}_\alpha(\overline{G})=2$. Also, $\lfloor 1/\alpha \rfloor \leq \lfloor 4/3 \rfloor=1$. Thus, $\mathsf{pd}_\alpha(G)+\mathsf{pd}_\alpha(\overline{G})=4=3+1 \leq \left\lceil \dfrac{4-1}{\lfloor \frac{1}{\alpha} \rfloor} \right\rceil +1$. Combining all these cases, we have the theorem.\qed

{\theorem For graphs $G$ and $\overline{G}$ without isolated vertices, $\mathsf{pd}_\alpha(G)+\mathsf{pd}_\alpha(\overline{G})\leq \left\lceil \dfrac{\lfloor n/2 \rfloor +2}{\lfloor \frac{1}{\alpha} \rfloor} \right\rceil +1$.}\\
\pf The theorem follows exactly as the proof of Theorem \ref{NG1} by using Proposition \ref{joseph-arumugam}.\qed

\section{$\alpha$-Partial Domination Spectrum of a Graph and its Consequences}

We define $\alpha$-partial domination spectrum of a graph $G$, denoted by $\mathsf{Sp}^p_\alpha(G)$, to be the set of distinct values of $\mathsf{pd}_\alpha(G)$ as $\alpha$ runs over $(0,1]$, i.e., $\mathsf{Sp}^p_\alpha(G)=\{\mathsf{pd}_\alpha: \alpha \in (0,1]\}$. Now, two cases may arise: either $\mathsf{Sp}^p_\alpha(G)$ is singleton or not. It is known that if for a graph $G$, $\gamma=1$, then $\mathsf{pd}_\alpha (G)= 1$ for all $\alpha \in (0,1]$, i.e., $|\mathsf{Sp}^p_\alpha(G)|=1$. On the other hand, if $\gamma \geq 2$, then $1,\gamma \in \mathsf{Sp}^p_\alpha(G)$, i.e., $|\mathsf{Sp}^p_\alpha(G)|\geq 2$.

Now, we move towards proving our main result that the $\alpha$-partial domination number changes its value only at rational points as $\alpha$ runs over $(0,1]$. However before doing that, we prove a lemma which we will use later. 
{\lemma \label{main-lemma} Let $G$ be a graph such that $|\mathsf{Sp}^p_\alpha(G)|>1$. Let $q\in \mathsf{Sp}^p_\alpha(G)$ such that $q\neq 1$. Let $A=\{\alpha \in (0,1]:\mathsf{pd}_\alpha(G)<q\}$ and $B=\{\alpha \in (0,1]:\mathsf{pd}_\alpha(G)\geq q\}$. Then there exists a rational number $\alpha^*\in (0,1)$ such that $A=(0,\alpha^*]$ and $B=(\alpha^*,1]$.}\\
\\
\pf   Since $q \neq 1$, $q$ is not the least element of $\mathsf{Sp}^p_\alpha(G)$. 
Observe that both $A$ and $B$ are non-empty, because $(0,\frac{\Delta+1}{n}]\subseteq A$ and $(1-\frac{1}{n},1]\subseteq B$. In fact, both $A$ and $B$ are intervals. It follows from the fact that $\alpha'<\alpha''$ implies $\mathsf{pd}_{\alpha'}(G)\leq \mathsf{pd}_{\alpha''}(G)$. Moreover, from the definition, it follows that $A\cup B=(0,1]$ and $A\cap B=\emptyset$. Thus there exists $\alpha^* \in (0,1]$ such that either $A=(0,\alpha^*), B=[\alpha^*,1]$ or $A=(0,\alpha^*], B=(\alpha^*,1]$. 

Claim 1: We claim that both $A$ and $B$ are left-open, right-closed intervals, i.e., $A=(0,\alpha^*], B=(\alpha^*,1]$. If possible, let $A=(0,\alpha^*), B=[\alpha^*,1]$. Let $p\in \mathsf{Sp}^p_\alpha(G)$ be the largest element in $\mathsf{Sp}^p_\alpha(G)$ less than $q$. Thus
there exists $\alpha' \in A$ such that $\mathsf{pd}_{\alpha'}=p$. This imply that for all $\alpha \in [\alpha',\alpha^*), \mathsf{pd}_\alpha=p$. Let $(\alpha_k)$ be a strictly monotonically increasing sequence in $[\alpha',\alpha^*)$ such that $(\alpha_k)$ converges to $\alpha^*$. Now as $\alpha_k \in [\alpha',\alpha^*)$, we have $\mathsf{pd}_{\alpha_k}=p$, i.e., for each $\alpha_k$, there exists $S_k \subseteq V$ with $|S_k|=p$ such that 
$|N[S_k]|\geq n\cdot \alpha_k.$ 
Moreover as $|S_k|=p<q$, $S_k$ is not a $\mathsf{pd}_{\alpha^*}$-set, i.e., $|N[S_k]|< n\cdot \alpha^*.$
Thus we get a sequence of $S_k$ of subsets of $V$ such that 
\begin{equation}\label{sandwich-equation}
\alpha_k \leq \dfrac{|N[S_k]|}{n}< \alpha^*, \forall k \in \mathbb{N}
\end{equation}
As $G$ is a finite graph, the number of choices for subsets $S_k$ of size $p$ is finite. Thus, the sequence $\left(\frac{|N[S_k]|}{n}\right)$ assumes finitely many values. Now, since $(\alpha_k)$ converges to $\alpha^*$, by Sandwich Theorem, the sequence $\left(\frac{|N[S_k]|}{n}\right)$ converges to $\alpha^*$. As any convergent sequence taking finitely many values is eventually constant, we have $\left(\frac{|N[S_k]|}{n}\right)$ to be an eventually constant sequence. Thus there exists $t \in \mathbb{N}$ such that $\frac{|N[S_k]|}{n}=\alpha^*$ for all $k \geq t$. This is a contradiction to Equation \ref{sandwich-equation}. Thus our claim is justified and hence $A=(0,\alpha^*]$ and $B=(\alpha^*,1]$.

Claim 2: $\alpha^* \in \mathbb{Q}$. If possible, let $\alpha^* \in (0,1]\setminus \mathbb{Q}$. We observe that $\mathsf{pd}_{\alpha^*}(G)=p$, because if $\mathsf{pd}_{\alpha^*}(G)<p$, then for all $\alpha \in A$, $\mathsf{pd}_\alpha(G)<p$ which contradicts the fact that $p\in \mathsf{Sp}^p_\alpha(G)$. Since, $\alpha^*$ is an irrational number, $n\alpha^* $ is not an integer. Now as $(0,1]\cap (\mathbb{R}\setminus \mathbb{Q})$ is dense in $(0,1]$, there exists an irrational number $\overline{\alpha}\in (0,1]\cap (\mathbb{R}\setminus \mathbb{Q})$ with $\overline{\alpha}>\alpha^*$ such that $\lceil n \alpha^* \rceil=\lceil n \overline{\alpha} \rceil$. (We omit the details of the proof)

Now let $S$ be a $\mathsf{pd}_{\alpha^*}$-set of $G$. Then $|N[S]|\geq n\alpha^*$. As $\alpha^* \in \mathbb{R}\setminus \mathbb{Q}$, we have $|N[S]|\geq \lceil n\alpha^* \rceil = \lceil n \overline{\alpha}\rceil\geq n\overline{\alpha} $. Thus $S$ is a $\overline{\alpha}$-partial dominating set in $G$ and $|S|=p$. 

On the other hand, as $\overline{\alpha}>\alpha^*,\overline{\alpha} \in B$. But $S$ being a $\overline{\alpha}$-partial dominating set of $G$ must have cardinality $\geq q$ (by definition of $B$). This is a contradiction. Hence $\alpha^* \in \mathbb{Q}$. \qed

Now we are in a position to prove the following theorem.
{\theorem \label{main-theorem} Let $G$ be a graph such that $\mathsf{Sp}^p_\alpha(G)=\{a_1,a_2,\ldots,a_t\}$ with $1=a_1<a_2< \ldots <a_t=\gamma$ and $t>1$. Then there exists $t-1$ rational numbers $\alpha_1< \alpha_2< \ldots <\alpha_{t-1}$ in $(0,1)\cap \mathbb{Q}$ such that 
\begin{enumerate}
\item $\forall \alpha \in (0,\alpha_1]$, $\mathsf{pd}_\alpha(G)=1$.
\item $\forall i \in \{1,2,\ldots,t-2\}, \forall \alpha \in (\alpha_i,\alpha_{i+1}]$, $\mathsf{pd}_\alpha(G)=a_{i+1}$.
\item $\forall \alpha \in (\alpha_{t-1},1]$, $\mathsf{pd}_\alpha(G)=\gamma$.
\end{enumerate}
}
\pf Substituting $q=a_t$ in Lemma \ref{main-lemma}, we get a rational number $\alpha_{t-1}$ such that $A_1=\{\alpha \in (0,1]:\mathsf{pd}_\alpha(G)< a_t\}=(0,\alpha_{t-1}]$ and $B_1=\{\alpha \in (0,1]:\mathsf{pd}_\alpha(G)\geq a_t\}=(\alpha_{t-1},1]$. However, as $a_t$ is the largest element in $\mathsf{Sp}_\alpha(G)$, we have $B_1=\{\alpha \in (0,1]:\mathsf{pd}_\alpha(G)= a_t\}=(\alpha_{t-1},1]$. 

Again, substituting $q=a_{t-1}$ in Lemma \ref{main-lemma}, we get a rational number $\alpha_{t-2}$ such that $A_2=\{\alpha \in (0,1]:\mathsf{pd}_\alpha(G)< a_{t-1}\}=(0,\alpha_{t-2}]$ and $B_2=\{\alpha \in (0,1]:\mathsf{pd}_\alpha(G)\geq a_{t-1}\}=(\alpha_{t-2},1]$. However, as $a_t$ and $a_{t-1}$ are the only two elements in $\mathsf{Sp}^p_\alpha(G)$ which are greater or equal to $a_{t-1}$ and $B_1=\{\alpha \in (0,1]:\mathsf{pd}_\alpha(G)= a_t\}=(\alpha_{t-1},1]$, we have $\{\alpha \in (0,1]:\mathsf{pd}_\alpha(G)= a_{t-1}\}=(\alpha_{t-2},\alpha_{t-1}]$. 

Continuing in this way, at one stage we substitute $q=a_2$ in Lemma \ref{main-lemma} to get a rational number $\alpha_1$ such that $A_{t-1}=\{\alpha \in (0,1]:\mathsf{pd}_\alpha(G)< a_2\}=(0,\alpha_{1}]$ and $B_{t-1}=\{\alpha \in (0,1]:\mathsf{pd}_\alpha(G) \geq a_{2}\}=(\alpha_1,1]$. By similar argument as that of above, we get $\{\alpha \in (0,1]:\mathsf{pd}_\alpha(G)= a_{2}\}=(\alpha_{1},\alpha_{2}]$. Moreover, as $a_1=1$ is the only value in $\mathsf{Sp}^p_\alpha(G)$ which is less than $a_2$, we have $A_{t-1}=\{\alpha \in (0,1]:\mathsf{pd}_\alpha(G)= a_1\}=(0,\alpha_{1}]$.
Hence the theorem.\qed

We call the $\alpha_i$'s obtained in Theorem \ref{main-theorem} as {\it critical values} of $\alpha$. Theorem \ref{main-theorem} has an immediate corollary.

{\corollary Let $G$ be a graph and $\overline{\alpha}$ be a irrational number in $(0,1)$. Then there exists $\epsilon >0$, such that for all $\alpha \in (\overline{\alpha}-\epsilon,\overline{\alpha}+\epsilon)$, $\mathsf{pd}_\alpha(G)$ is constant.}\\
\\ \pf The corollary follows from Theorem \ref{main-theorem} and denseness of rationals and irrationals in $\mathbb{R}$. \qed

Our next goal is to find an upper bound on the size of the $\alpha$-partial domination spectrum of a graph. Before that we prove a lemma.

{\lemma \label{bound-attain-lemma} Let $G$ be a graph such that $\mathsf{Sp}^p_\alpha(G)=\{a_1,a_2,\ldots,a_t\}$ with $a_1<a_2< \ldots <a_t$ and let $\alpha_i$'s be as in Theorem \ref{main-theorem}. Then for each $\alpha_i$, there exists a $\mathsf{pd}_{\alpha_i}$-set $S_i \subseteq V$ such that $$\alpha_i=\frac{|N[S_i]|}{n}.$$}
\pf Since $S_i$ is a $\mathsf{pd}_{a_i}$-set of $G$, we have $|S_i|=a_i$ and 
\begin{equation}\label{attaining-equation}
\frac{|N[S_i]|}{n}\geq \alpha_i.
\end{equation}
If possible, the inequality in Equation \ref{attaining-equation} is strict. But in that case, by denseness of real numbers, we can find $\alpha'>\alpha_i$ such that $\frac{|N[S_i]|}{n}\geq \alpha'> \alpha_i$. Thus $S_i$ is $\alpha'$-partial dominating set of $G$ and hence $\mathsf{pd}_{\alpha'}(G)\leq |S_i|=a_i<a_{i+1}$. However as $\alpha'> \alpha_i$, we have $\mathsf{pd}_{\alpha'}(G)\geq a_{i+1}$. This is a contradiction. Thus, there exists $S_i$ such that Equation \ref{attaining-equation} holds with equality. Hence the theorem. \qed

{\theorem \label{spectrum-upper-bound} For any graph $G$ without isolated points, the critical values belong to the set $\{\frac{\Delta +1}{n},\frac{\Delta +2}{n},\cdots, \frac{n-1}{n}\}$ and $|\mathsf{Sp}^p_\alpha(G)|\leq n-\Delta$. }\\
\\ \pf By Lemma \ref{bound-attain-lemma}, for every critical value $\alpha_i$, there exists a $\mathsf{pd}_{\alpha_i}$-set $S_i \subseteq V$ such that $$\alpha_i=\frac{|N[S_i]|}{n}.$$
Thus the first part of the theorem follows from Proposition \ref{basic-bound-1}, \ref{basic-bound-2} and the observation that $|N[S_i]|\leq n$.

For the second part, observe that $|\mathsf{Sp}^p_\alpha(G)|$ is one more than the number of critical values. Thus, $|\mathsf{Sp}^p_\alpha(G)|\leq 1+(n-\Delta-1)= n-\Delta$.\qed

{\remark The upper bound given in Theorem \ref{spectrum-upper-bound} is tight: Consider an $n$ vertex graph which consists of a clique and some isolated vertices.}

\section{Conclusion}
In this paper, we introduced a new graph invariant called the partial domination number of a graph. From an applications standpoint, it  can be interpreted as the measure of the maximum surveillance possible if a fraction of minimum number of guards needed is available.  We studied different bounds on the partial domination number of a graph $G$ with respect to its order, maximum degree, domination number etc. 

\section*{Acknowledgement}
The author is grateful to Geertrui Van de Voorde from University of Ghent for pointing out a mistake in an earlier version of the paper. We note that, authors in \cite{laskar-partial-domination} independently, has proposed the same notion of partial domination. However, their research is mainly focussed on the case $\alpha=\frac{1}{2}$, whereas we dealt with general values of $\alpha$.
The research is partially funded by NBHM Research Project Grant, (Sanction No. 2/48(10)/2013/ NBHM(R.P.)/R\&D II/695), Government of India.

\end{document}